\documentclass[11pt]{article}
\usepackage[utf8]{inputenc}
\usepackage{color}
\usepackage{xcolor}

\usepackage{amsmath,amsfonts,euscript,amssymb,amsthm,a4,times}

\allowdisplaybreaks

\def\XXint#1#2#3{{\setbox0=\hbox{$#1{#2#3}{\int}$}
\vcenter{\hbox{$#2#3$}}\kern-.5\wd0}}

\newtheorem{theorem}{Theorem}[section]

\newtheorem{lemma}[theorem]{Lemma}

\theoremstyle{definition}
\newtheorem{definition}[theorem]{Definition}

\makeatletter \catcode`@=11
\newbox\tr@tto
\setbox\tr@tto=\hbox{{\count0=0\dimen0=-,9pt\dimen1=1,1pt\loop\ifnum
    \count0<11 \advance \count0 by1 \vrule width.51pt height\dimen1
    depth\dimen0\kern-0.17pt\advance\dimen0 by-0.05pt\advance\dimen1
    by0.1pt\repeat \loop\ifnum\count0<21\advance \count0 by1 \vrule
    width.6pt height\dimen1 depth\dimen0\kern-0.2pt \advance\dimen0
    by-0.1pt\advance\dimen1 by 0.05pt\repeat}}
\def\medint{\displaystyle\copy\tr@tto\kern-10.4pt\int}
\numberwithin{equation}{section}

\def\R{{\mathbb R}}
\def\RN{{\mathbb R}^{N}}
\def\M{{\mathbb R}^{N\times n}}

\def\proofof#1{\begin{proof}[Proof of #1]}

\catcode`@=11
\newbox\tr@tto
\setbox\tr@tto=\hbox{{\count0=0\dimen0=-,9pt\dimen1=1,1pt\loop\ifnum\count0<11 \advance \count0 by1 \vrule width.51pt height\dimen1
     depth\dimen0\kern-0.17pt\advance\dimen0 by-0.05pt\advance\dimen1
     by0.1pt\repeat \loop\ifnum\count0<21\advance \count0 by1 \vrule
     width.6pt height\dimen1 depth\dimen0\kern-0.2pt \advance\dimen0
     by-0.1pt\advance\dimen1 by 0.05pt\repeat}}
\def\medint{\displaystyle\copy\tr@tto\kern-10.4pt\int}
\catcode`@=12

\newcommand{\LL}{\mathrm{L}}
\newcommand{\WW}{\mathrm{W}}

\numberwithin{equation}{section}

\begin{document}
\title{Lipschitz regularity of minimizers of
 variational integrals with variable exponents}
\author{Michela Eleuteri- Antonia Passarelli di Napoli}

\maketitle

\begin{abstract}
In this paper we prove the Lipschitz regularity  for local minimizers of  convex
variational integrals of the form
\[
\mathfrak{F}( v, \Omega )= \int_{\Omega} \! F(x, Dv(x)) \, dx,
\]
where, for ${n > 2}$ and $N\ge 1$, $\Omega$ is a bounded open set in $\mathbb{R}^n$, {$u \in \WW^{1,1}(\Omega , \RN )$ and the energy density $F:\Omega\times \M\to \mathbb{R}$  satisfies  the so called variable growth conditions.} The main novelty of the paper is that we assume an almost critical regularity in the Orlicz Sobolev setting for the energy density as a function of the $x$ variable.
\end{abstract}

\maketitle

\noindent
{\footnotesize {\bf AMS Classifications.}  49N15; 49N60; 49N99.}

\noindent
{\footnotesize {\bf Key words and phrases.}  Variational integrals; Variable exponents; Regularity of minimizers.}

\bigskip

\section{Introduction}
\bigskip
The aim of this paper is to prove a Lipschitz regularity result for the local minimizers of  convex
variational integrals of the form
\begin{equation}\label{defint}
\mathfrak{F}( v, \Omega )= \int_{\Omega} \! F(x, Dv(x)) \, dx,
\end{equation}
where, for ${n > 2}$ and $N\ge 1$, $\Omega$ is a bounded open set in $\mathbb{R}^n$, {$u \in \WW^{1,1}(\Omega , \RN )$ and $F:\Omega\times \M\to \mathbb{R}$  is a strictly convex $C^2$ function with respect to the second variable satisfying  the so called variable growth conditions.}
 More precisely, we shall assume that there exist functions $p:\Omega\to (1,+\infty)$, $k:\Omega\to [0,+\infty)$  and  constants $\ell, \nu>0$, $0 < C_1 < C_2$ such that
$$C_1 \, (1+|\xi|^2)^{\frac{p(x)}{2}}\le F(x,\xi)\le  \, C_2 \, (1+|\xi |^2)^{\frac{p(x)}{2}} \eqno{\rm (H1)}$$
$$\langle {F_{\xi \xi}} (x,\xi)\eta,\eta\rangle \ge \nu  (1 + |\xi |^2)^{\frac{p(x)-2}{2}}|\eta|^2 \eqno{\rm (H2)}$$
$$|{F_{\xi \xi}}(x,\xi)| \le  \ell (1 + |\xi |^2)^{\frac{p(x)-2}{2}}
\eqno{\rm (H3)}$$
$$|{F_{x \xi}}(x,\xi)|\le \, k(x)(1+|\xi |^2)^{\frac{p(x)-1}{2}}{\log (e+|\xi|^2)}\,,\eqno{\rm (H4)}$$
for all $\xi,\eta \in \M$ and for a.e. $x\in \Omega$.
\\
For the functions $p, k$ appearing in previous assumptions, we   will assume that  there exists $\theta>n$ such that
$$
    k\in L^n\log^{n+\theta}_{\mathrm{loc}}(\Omega) \eqno{\rm (H5)}$$
and moreover
   $$ p(x)\in W^{1,1}_{\mathrm{loc}}(\Omega)\,\, \text{with}\,\, |Dp|\in L^n\log^{n+\theta}_{\mathrm{loc}}(\Omega).\eqno{\rm (H6)}$$
It is worth pointing out that assumption (H6) implies that $p(x)$ is a  continuous function with a logarithmic modulus of continuity, by virtue of the Sobolev embedding theorem in Orlicz space (for more details we refer to Sections \ref{Orlicz} and \ref{Orlicz2}). Such modulus of continuity for the exponent $p$ cannot be removed as well as, by virtue of the counterexample by Zhikov \cite{Z5} the Lavrentiev phenomenon may appear and therefore regularity of the minimizers cannot be in general expected.
\\
It is worth mentioning that the continuity of the exponent $p(x)$, with a logarithmic modulus of continuity, has been revealed sufficient to prove the H\"older continuity of the local $\mathfrak{F}$-minimizers  (see \cite{AM1, cos-min, ele1, GPT, GPTR}), while the regularity of their gradient usually has been obtained imposing the H\"older continuity of $p(x)$. As far as we know, our main Theorem  is the first regularity result for the gradient of the local $\mathfrak{F}$-minimizers obtained assuming a Sobolev regularity on $p(x)$ that doesn't imply its H\"older continuity.
\\
In what follows, in order to have that the model case $$F(x,\xi)=(1+|\xi|^2)^{\frac{p(x)}{2}}$$ 
satisfies our assumptions, we shall suppose, without loss of generality, that
\begin{equation}
   |Dp(x)|\le k(x)\qquad\qquad \text{a.e. in}\,\,\Omega.
\end{equation}
We also remark that functions of the form 
$$F(x,\xi)=a(x)(1+|\xi|^2)^{\frac{p(x)}{2}},$$ 
with $a(x)$ bounded and weakly differentiable and such that $Da \in L^n\log^{\theta+n}L_{\mathrm{loc}}(\Omega)$, satisfy (H1)--(H6) with $$k(x)=\max\{|Da(x)|, |Dp(x)|\}.$$
Finally, in order to avoid the irregularity phenomena peculiar of the vectorial case, we suppose that there exists $g:\Omega\rightarrow  [0,+\infty)\to [0,+\infty)$ such that
\begin{equation}\label{modulo}F(x,\xi)=g(x,|\xi|).
\end{equation}
Lipschitz continuity of local minimizers with convex integrands with standard $p$-growth is a classic topic and for an exhaustive treatment we refer to \cite{giaquinta,Giusti} and the references therein.
Many contributions are available also in case of autonomous integrands with non standard $(p,q)$ growth conditions since the pioneering papers \cite{mar91, mar93, mar96} by Marcellini  appeared, see also the more recent results in \cite{M21}, \cite{M21bis}.
\\
When dealing with non standard growth conditions in the non autonomous case, i.e.  when the energy densities depend also on the $x$-variable, things become more involved since the Lavrentiev phenomenon may appear.
\\
However, in the last few years, many Lipschitz regularity results have been established in case the partial map $x \mapsto F(x,\xi)$ belongs to the Sobolev class $W^{1,r}_{\mathrm{loc}}(\Omega)$, with $r>n$ and therefore, by the Sobolev imbedding Theorem, is H\"older continuous (see \cite{EMM1,EMM2,DeFilippis-Mingione 2020}).
\\
Here we deal with the $p(x)$-growth conditions, that, in some sense, is intermediate between the standard and the non standard growth cases and which is nowadays a very popular topic in the Calculus of Variations and Nonlinear Analisis (for more details see \cite{DHHR, Harjulehto} and the references therein).
\\
Our aim is to show that the Lipschitz continuity of the local minimizers  holds true under assumption (H5) which is weaker than that made in \cite{EMM2} and analogous to that made in \cite{CGHP20} for the standard growth case. Actually it is well known that
\[
W^{1,r} \subset W^{1, L^n \log^{\beta} L} \subset W^{1,n} \qquad \forall r > n, \,\, \forall \beta > 0.
\]
It worth pointing out that assumption (H5) has been already employed in \cite{GP, FG} to obtain second order regularity of minimizers.
\\
Before stating our main result, we recall the definition of local minimizer
\begin{definition}
A mapping $u\in \WW^{1,1}_{\rm loc}(\Omega, \RN )$ is a local $\mathfrak{F}$--minimizer if
$F(x,Du) \in \LL^{1}_{\rm loc}(\Omega )$ and
$$
\int_{\mathrm{supp}\varphi} \! F(x,Du) \, dx \leq \int_{\mathrm{supp}\varphi} \! F(x,Du+D\varphi) \, dx
$$
for  any $\varphi\in {C}_0^{\infty}(\Omega, \RN )$.
\end{definition}

{The main result of the paper can be stated as follows.
\begin{theorem}
\label{main}
Let $u$ be a local minimizer of the functional \eqref{defint} under
the assumptions {\rm {(H2)--(H6)}}. Suppose moreover that there exists a constant $M_1 > 0$ such that
\begin{equation}
\label{M1}
\int_{\mathcal{O}} |Du|^{p(x)} \, dx \le \, M_1
\end{equation}
{for all open $\mathcal{O} \subset \Omega.$}
Then $u$ is locally Lipschitz continuous and the following estimate holds for any {$B_{R_0} \Subset \mathcal{O}$}
\[
{||1+Du||_{L^\infty \left(B_{\frac{R_0}{2}}\right)}\le  C\left (\int_{B_{R_0}}  \, (1 + |Du|^2)^{\frac{p^-}{2}} dx \right) ^{\frac{2}{p^-}},}
\]
{with $C \equiv C (n,\theta, \ell, \nu, M_1, {||k||_{L^n\log^{n+\theta}L(B_{R_0})}}, R_0, p^-,p^+),$ where $p^-=\inf_\Omega p(x)$, $p^+=\sup_\Omega p(x)$ .}
\end{theorem}
}
The proof relies, as usual, on the combination of a suitable a priori estimate with an approximation argument. The a priori estimate is established through the well known Moser iteration argument. The novelty here is to show that such argument works under our weaker assumption on the partial maps $x \mapsto F(x, \xi),$ $x \mapsto F_{\xi}(x, \xi)$. To this  aim, we will  plenty use the properties of Orlicz-Zygmund classes and variable Lebesgue spaces, as well as the Sobolev embedding theorem in the limit case.
\\
Once the a priori estimate is established, the approximation can be made effective arguing as done in the $(p,q)$ growth case.
\\
We conclude mentioning that the arguments of the proof of Theorem \ref{main} can be easily adapted to prove higher differentiability of the local $\mathfrak{F}$--minimizers, following the lines of \cite{CEP}.
\\
The paper is organized as follows: Section \ref{due} contains some notations and preliminary results; Section \ref{tre} is devoted to the proof of Theorem \ref{main} and finally Section \ref{quattro} is concerned with the approximation and passage to the limit.

\section{Notations and preliminary results}

\label{due}

{In this paper we shall denote by $C$  a
general positive constant that may vary on different occasions, even within the
same line of estimates.
Relevant dependencies  will be suitably emphasized using
parentheses or subscripts.  In what follows, $B(x,r)=B_r(x)=\{y\in \R^n:\,\, |y-x|<r\}$ will denote the ball centered at $x$ of radius $r$.
We shall omit the dependence on the center and on the radius when no confusion arises.
\\
We recall the following Lemma, whose proof can be found for example in \cite[Lemma 6.1]{Giusti}.
 \begin{lemma}\label{holf} Let $h:[\rho_0, R_{0}]\to \mathbb{R}$ be a non-negative bounded function and $0<\vartheta<1$,
$A, B\ge 0$ and $\beta>0$. Assume that
$$
h(s)\leq \vartheta h(t)+\frac{A}{(t-s)^{\beta}}+B,
$$
for all $\rho_0\leq s<t\leq R_{0}$. Then
$$
h(r)\leq C \, \frac{A}{(R_{0}-\rho_0)^{\beta}}+ C \, B ,
$$
where $C=C(\vartheta, \beta)>0$.
\end{lemma}}

\subsection{The Orlicz-Zygmund space $L \log^{\beta} L$}
\label{Orlicz} The Orlicz-Zygmund space $L \log^{\beta} L(\Omega; \mathbb{R}^N)$ is defined as
\[
L \log^{\beta} L(\Omega; \mathbb{R}^N) := \left \{f \in L^1(\Omega; \mathbb{R}^N): \, \int_{\Omega} |f| \log^{\beta} (e + |f|) \, dx < \infty\right \}
\]
and it becomes a Banch spaces with the Luxemburg norm
\[
\|f\|_{L \log^{\beta} L(\Omega)} := \inf \left \{ \lambda > 0: \medint_{\Omega} \left |\frac{f}{\lambda} \right | \log^{\beta} \left (e + \left | \frac{f}{\lambda} \right | \right ) \, dx \le \, 1\right \}
\]
For more details on Orlicz spaces we refer to \cite{RR91} and for further properties of the Orlicz-Zygmund classes we refer to \cite{IMMP}.
\\
In particular, this space embeds in any $L^p(\Omega; \mathbb{R}^N)$ for $p > 1$; more precisely for any $p > 1$ the following inequality takes place
\begin{equation}\label{embed}
\|f\|_{L \log^{\beta} L(\Omega)} \le \, C \left ( \medint_{\Omega} |f|^p \, dx \right )^{1/p} \qquad \forall f \in L^p (\Omega; \mathbb{R}^N)
\end{equation}
where the constant $C$ depends only on $p$ and blows up when $p \rightarrow 1$.
\\
By a characherization due to Iwaniec \cite{I92}, if we set
\[
[f]_{L \log^{\beta} L(\Omega)} := \medint_{\Omega} |f| \log^{\beta} \left (e + \frac{|f|}{\|f\|_1} \right ) \, dx
\]
where
\[
\|f\|_1 := \medint_{\Omega}|f| \,dx
\]
then it turns out that the quantity 
$[f]_{L \log^{\beta} L(\Omega)}$, for $\beta\ge 0$, is comparable with the Luxemburg norm in $L \log^{\beta} L(\Omega; \mathbb{R}^N),$ in the sense that there exists a constant $C \equiv C(\beta) \ge \, 1$, independent of $\Omega$ and $f$, such that
\[
C^{-1} \|f\|_{L \log^{\beta} L(\Omega)} \le \, [f]_{L \log^{\beta} L(\Omega)} \le \, C \, \|f\|_{L \log^{\beta} L(\Omega)}
\]
Therefore, in view of \eqref{embed}, we get 
\begin{equation}
\label{lemmaAM}
\medint_{\Omega} |f| \log^{\beta} \left (e + \frac{|f|}{\|f\|_1} \right ) \, dx \le \, C(p, \beta) \left (\medint_{\Omega} |f|^p \, dx \right )^{1/p} \qquad \forall f \in L^p (\Omega; \mathbb{R}^N)
\end{equation}
for every  $\beta \ge 0$.


Similarly, one can define $L^p\log^\beta L$ as the space of functions such that
\[
[f]_{L^p \log^{\beta} L(\Omega)} := \medint_{\Omega} |f|^p \log^{\beta} \left (e + \frac{|f|}{\|f\|_p} \right ) \, dx<+\infty
\]
where
\[
\|f\|_p := \left(\medint_{\Omega}|f|^p \,dx\right)^{\frac{1}{p}}
\]
whose properties and embeddings can be deduced as done before for the case $p=1$.

We conclude this section with the following elementary inequality which will be useful in the sequel
\begin{lemma}\label{elem}
For every $s,t>0$, and for every $\beta>0$ we have
$$st\le \delta s\log^\beta \left(e+\tau s\right)+\frac{t}{\tau}\left(\exp\left(\frac{t}{\delta}\right)^{\frac{1}{\beta}}-1\right) ,$$
where $\delta,\tau>0.$
\end{lemma}

\subsection{Sobolev embeddings in Orlicz spces}
\label{Orlicz2}

We recall the following  
embedding theorem whose proof can be found in \cite{Cianchi}.
%

\begin{theorem}\label{sob2}
Let $h\in W^{1,1}_0(\Omega)$ be a function such that $|Dh|\in L^n\log^\sigma L(\Omega)$, some $\sigma >n-1$. Then $h\in C^0(\Omega)$ and 
\begin{equation}
   |h(x)- h(y)|\le \frac{C(n)}{\log^{\frac{\sigma-n+1}{n}}\left(e+\frac{1}{|x-y|}\right)} ||Dh||_{L^n\log^\sigma L(\Omega)},
\end{equation}
for some constant $C \equiv C(n).$
\end{theorem}   
 
\subsection{Spaces with variable exponents}

Let $p: \Omega\to (1,+\infty)$ be the function appearing in assumptions (H1)--(H4). By virtue of assumption (H6) and Theorem \ref{sob2}, we have that $p$ is continuous with a logarithmic modulus of continuity, i.e.
$$ |p(x)- p(y)|\le \frac{C(n)}{\log^{\frac{\theta+1}{n}}\left(e+\frac{1}{|x-y|}\right)} ||Dp||_{L^{n}\log^{n+\theta} L_{\mathrm{loc}}(\Omega)}. $$ Due to the local nature of the results we want
to prove, it is not restrictive to assume the global boundedness of the function $p(x)$, i.e. that
there exist constants $p^-, p^+$ such that
\begin{equation}\label{p+p-}
    1<p^- = \inf_{\Omega} p(x) \le p(x)\le p^+ = \sup_{\Omega} p(x) <+\infty\,.
\end{equation}

In the sequel we need to deal {with} the Lebesgue spaces with variable exponents, $L^{p(\cdot)}(\Omega; \mathbb{R}^N);$  for more details we refer to the monographs \cite{DHHR, CUF}.



\begin{definition}The space $L^{p(x)}(\Omega; \mathbb{R}^N)$ is defined as
\begin{equation}
    L^{p(x)}(\Omega; \mathbb{R}^N)=\left\{v:\Omega\to \mathbb{R}^n\,:\, v\,\, \text{is measurable and}\,\, \int_{\Omega}|v|^{p(x)}\,dx<+\infty\right\}
\end{equation}
and, once equipped with the Luxemburg norm, it becomes a Banach space.

\end{definition}   

Next, we recall the definition of Sobolev space with variable exponents.
\begin{definition}The space $W^{1,p(x)}(\Omega; \mathbb{R}^N)$
is defined as
$$W^{1,p(x)}(\Omega; \mathbb{R}^N):=\{v\in L^{p(x)}(\Omega; \mathbb{R}^N)\,\,:\,\, Dv\in L^{p(x)}(\Omega; \mathbb{R}^{nN})\}$$
where $Dv$ denotes the distributional gradient of $v$.
\end{definition}   

It is clear from (H1) that 
$$\int_{\mathcal{O}}F(x,Du)\,dx<+\infty, \quad \forall \mathcal{O}\Subset\Omega\,\,\Longleftrightarrow u\in W^{1,p(x)}_{\mathrm{loc}}(\Omega; \mathbb{R}^N) .$$



\section{Proof of Theorem \ref{main}}

\label{tre}

{The proof of the main result will be established through two main steps: the a priori estimate and the approximation procedure, together with a localization process, aimed to apply some known results in the non-standard growth setting.}

\subsection{The a priori estimate}
{As usual, we are going to establish an a priori estimate for local minimizer $u$ of $\mathfrak{F}$ such that $u\in W^{2,2}_{\mathrm{loc}}(\Omega)\cap W^{1,\infty}_{\mathrm{loc}}(\Omega)$.} We start from the second variation of the functional $\mathfrak{F}$, i.e., for a fixed $s\in\{1,\dots,n\}$
\begin{equation}
\label{second-variation}
\int_\Omega\left(\sum_{i,j,\alpha,\beta}F_{\xi_i^\alpha\xi_j^\beta}(x,Du)\varphi^\alpha_{x_i}u^\beta_{x_jx_s}+\sum_{i,\alpha}F_{\xi_i^\alpha x_s}(x,Du)\varphi^\alpha_{x_i}\right)\,dx=0.
\end{equation}
{Fix a ball $B_{R_0}\Subset \Omega$ and choose radii $0 < \frac{R_0}{2} < \rho < s < t < R < R_0$; let $\eta\in C^1(B_{t})$ be a cut-off function such that $0 \le \eta \le 1$, $\eta = 1$ on $B_s$ and $D \eta \le \, \frac{C}{t - s}.$ For}  $\gamma>0$,  we choose $$\varphi^{\alpha}=\eta^2u^\alpha_{x_s}(1+|Du|^2)^{\frac{\gamma}{2}}$$
as test function in previous equality, {which is legitimate by virtue of the a priori assumption on $u$}. We compute
\begin{eqnarray*}
\varphi_{x_i}^{\alpha} &=& 2 \eta \eta_{x_i} u_{x_s}^{\alpha} (1 + |Du|^2)^{\frac{\gamma}{2}} + \eta^2 u_{x_s x_i}^{\alpha} (1 + |Du|^2)^{\frac{\gamma}{2}} \\
&& + \eta^2 u_{x_s}^{\alpha} \gamma (1 + |Du|^2)^{\frac{\gamma}{2} - 1} |Du| (|Du|)_{x_i}.
\end{eqnarray*}
Plugging this expression in \eqref{second-variation}, we obtain
\begin{eqnarray}
0 &=& \int_{\Omega} 2 \eta (1 + |Du|^2)^{\frac{\gamma}{2}} \sum_{i,j,\alpha, \beta} F_{\xi_i^{\alpha} \xi_j^{\beta}}(x, Du) \eta_{x_i} u_{x_s}^{\alpha} u_{x_s x_j}^{\beta} \, dx \nonumber\\
&& + \int_{\Omega} \eta^2 (1 + |Du|^2)^{\frac{\gamma}{2}} \sum_{i,j,\alpha, \beta} F_{\xi_i^{\alpha} \xi_j^{\beta}}(x, Du) u_{x_s x_i}^{\alpha} u_{x_s x_j}^{\beta} \, dx \nonumber\\
&& + \int_{\Omega} \eta^2 \, \gamma \, (1 + |Du|^2)^{\frac{\gamma}{2} - 1} \sum_{i,j,\alpha, \beta} F_{\xi_i^{\alpha} \xi_j^{\beta}}(x, Du) u_{x_s}^{\alpha} u_{x_s x_j}^{\beta} |Du| \, (|Du|)_{x_i}  \, dx\nonumber\\
&& + \int_{\Omega} 2 \eta (1 + |Du|^2)^{\frac{\gamma}{2}} \sum_{i, \alpha} F_{\xi_i^{\alpha} x_s}(x, Du)\eta_{x_i} u_{x_s}^{\alpha} \, dx\nonumber\\
&& + \int_{\Omega} \eta^2 (1 + |Du|^2)^{\frac{\gamma}{2}} \sum_{i, \alpha} F_{\xi_i^{\alpha} x_s}(x, Du)  u_{x_s x_i}^{\alpha} \, dx \nonumber\\
&& + \int_{\Omega} \eta^2 \gamma(1 + |Du|^2)^{\frac{\gamma}{2} -1} \sum_{i, \alpha} F_{\xi_i^{\alpha} x_s}(x, Du) u_{x_s}^{\alpha} \, |Du| \, (|Du|)_{x_i}\, dx \nonumber\\
&=:& I_{1,s} + I_{2,s} + I_{3,s} + I_{4,s} + I_{5,s} + I_{6,s}. \label{I1-I6}
\end{eqnarray}

Next, we sum  all terms in the previous equation with respect to $s$ from 1 to $n$, and we denote by $I_1-I_6$ the corresponding integrals.

This yields
\begin{equation}
\label{somma1}
I_2 + I_3 \le \, |I_1| + |I_4| + |I_5| + |I_6|.
\end{equation}
We estimate $I_3$ by using assumption \eqref{modulo}. We first observe that
\[
F_{\xi_i^{\alpha} \xi_j^{\beta}}(x, \xi) = \left (\frac{g_{tt}(x, |\xi|)}{|\xi|^2} - \frac{g_t(x, |\xi|)}{|\xi|^3} \right )\xi_i^{\alpha} \xi_j^{\beta} + \frac{g_t(x, |\xi|)}{|\xi|} \delta_{\xi_i^{\alpha} \xi_j^{\beta}}.
\]
Therefore
\begin{eqnarray}
&& \sum_{i,j,s,\alpha, \beta} F_{\xi_i^{\alpha} \xi_j^{\beta}}(x, Du) u_{x_s}^{\alpha} u_{x_s x_j}^{\beta} (|Du|)_{x_i} \nonumber\\
&=& \left (	\frac{g_{tt}(x, |Du|)}{|Du|^2} - \frac{g_t(x, |Du|)}{|Du|^3} \right ) \sum_{i,j,s,\alpha,\beta} u_{x_s}^{\alpha} u_{x_s x_j}^{\beta} u_{x_j}^{\beta} u_{x_i}^{\alpha} (|Du|)_{x_i} \nonumber\\
&& + \frac{g_t(x, |Du|)}{|Du|} \sum_{s,i,\alpha}  u_{x_s}^{\alpha} u_{x_s x_i}^{\alpha} (|Du|)_{x_i}\label{I4pp}\\
&=& \left (	\frac{g_{tt}(x, |Du|)}{|Du|} - \frac{g_t(x, |Du|)}{|Du|^2} \right ) \sum_{\alpha} \left [\sum_i u_{x_i}^{\alpha} (|Du|)_{x_i} \right ]^2 + g_t(x, |Du|) |D(|Du|)|^2,\nonumber
\end{eqnarray}
where we used that
\[
(|Du|)_{x_i} = \frac{1}{|Du|} \sum_{\alpha, s} u_{x_i x_s}^{\alpha} u_{x_s}^{\alpha}.
\]
Thus, coming back to the estimate of $I_3$ from \eqref{I4pp} we deduce
\begin{eqnarray*}
I_3 &=& \int_{\Omega} \eta^2 (1 + |Du|^2)^{\frac{\gamma}{2} - 1} |Du|\\
&&\quad\cdot \Bigg \{ 
\left (	\frac{g_{tt}(x, |Du|)}{|Du|} - \frac{g_t(x, |Du|)}{|Du|^2} \right ) \sum_{\alpha} \left [\sum_i u_{x_i}^{\alpha} (|Du|)_{x_i} \right ]^2 \\
&& \qquad+ g_t(x, |Du|) |D(|Du|)|^2
\Bigg \} \, dx.
\end{eqnarray*}
As usual, by virtue of  the inequality 
$
|D(|Du|)|^2 \le \, |D^2 u|^2,
$
one can conclude that
\[
I_3 \ge \, \int_{\Omega} \eta^2  (1 + |Du|^2)^{\frac{\gamma}{2} - 1} |Du| \frac{g_{tt}(x, |Du|)}{|Du|} \sum_{\alpha} \left (\sum_i u_{x_i}^{\alpha} (|Du|)_{x_i}\right )^2 \, dx \ge 0,
\]
where we also used that $g_{tt}(x, |Du|) \ge 0$.
Therefore, estimate \eqref{somma1} implies
\begin{eqnarray}\label{ristart}
I_2\le |{I}_1|+|{I}_4|+|{I}_5|+|{I}_6|.
\end{eqnarray}
By the Cauchy-Schwarz inequality, the Young inequality and  assumption (H2), we have
\begin{eqnarray}\label{I1}
|{I}_1| &=& 2\left |\int_{\Omega}  \eta  (1 + |Du|^2)^{\frac{\gamma}{2}}   \sum_{i,j,s, \alpha, \beta} F_{\xi_i^{\alpha} \xi_j^{\beta}}(x, Du) u_{x_j x_s}^{\beta} \eta_{x_i} u_{x_s}^{\alpha} \, dx\right |\cr\cr
&\le& 2\int_{\Omega}  \eta  (1 + |Du|^2)^{\frac{\gamma}{2}} \cr\cr
&& \times \left \{ \sum_{i,j,s, \alpha, \beta} F_{\xi_i^{\alpha} \xi_j^{\beta}}(x, Du) \eta_{x_i} \eta_{x_j} u_{x_s}^{\alpha} u_{x_s}^{\beta} \right \}^{1/2} \cr\cr
&&\times\left \{ \sum_{i,j,s, \alpha, \beta} F_{\xi_i^{\alpha} \xi_j^{\beta}}(x, Du) u^{\alpha}_{x_s x_i} \, u^{\beta}_{x_s x_j}\right \}^{1/2} \, dx \cr\cr
&\le &  2   \int_{\Omega}  (1 + |Du|^2)^{\frac{\gamma}{2}}  \sum_{i,j,s, \alpha, \beta} F_{\xi^{\alpha}_i \xi^{\beta}_j}(x, Du) \eta_{x_i} \eta_{x_j} u_{x_s}^{\alpha} u_{x_s}^{\beta} \, dx\cr\cr
&&  + \frac{1}{2} \int_{\Omega} \eta^2  (1 + |Du|^2)^{\frac{\gamma}{2}}  \sum_{i,j,s, \alpha, \beta} F_{\xi_i^{\alpha} \xi_j^{\beta}}(x, Du) u^{\alpha}_{x_s x_i} \, u^{\beta}_{x_s x_j}  \, dx \cr\cr
&\le& 2\ell \int_{\Omega} |D \eta|^2 \,  (1 + |Du|^2)^{\frac{\gamma + p(x)}{2}}  \, dx \cr\cr
&& + \frac{1}{2} \int_{\Omega} \eta^2  (1 + |Du|^2)^{\frac{\gamma}{2}}  \, \sum_{i,j,s, \alpha, \beta} F_{\xi_i^{\alpha} \xi_j^{\beta}}(x, Du) u_{x_j x_s}^{\alpha} u_{x_i x_s}^{\beta} \, dx.
\end{eqnarray}

Using assumption (H4) {and observing that $\log(e + |Du|^2) \le \, 2 \log (e + |Du|)$,} we get
\begin{eqnarray*}
|I_4| &\le& 2 \, \int_{\Omega}  \eta  k(x)(1 + |Du|^2)^{\frac{\gamma+p(x)-1}{2}}\log(e+|Du|)|D\eta||Du| \, dx \cr\cr
&\le& 2 \,\int_{\Omega}  \eta  k(x)(1 + |Du|^2)^{\frac{\gamma+p(x)}{2}}\log(e+|Du|)|D\eta| \, dx\cr\cr
&\le &  \,\int_{\Omega}  \eta^2  k^2(x)(1 + |Du|^2)^{\frac{\gamma+p(x)}{2}}\log^2(e+|Du|) \, dx\cr\cr
&&\quad + \,  \, \int_{\Omega} |D \eta|^2  (1 + |Du|^2)^{\frac{\gamma+p(x)}{2}} \,dx,
\end{eqnarray*}
where, in the last line, we used Young's inequality.
Again by (H4) and Young's inequality we get, {for a suitable value of the parameter $\sigma > 0$ which will be determined later}
\begin{eqnarray*}
|I_5| + |I_6| &\le& 2 \, \int_{\Omega}  \eta^2  k(x)(1 + |Du|^2)^{\frac{\gamma+p(x)-1}{2}}\log(e+|Du|)|D^2u| \, dx \cr\cr
&\le& \sigma\int_{\Omega}  \eta^2  (1 + |Du|^2)^{\frac{\gamma+p(x)-2}{2}}|D^2u|^2 \, dx\cr\cr
& &\quad + \frac{2}{\sigma} \int_{\Omega}  \eta^2  k^2(x)(1 + |Du|^2)^{\frac{\gamma+p(x)}{2}}\log^2(e+|Du|) \, dx.
\end{eqnarray*}
Therefore
\begin{eqnarray}
|I_4| + |I_5| + |I_6| &\le&  \sigma\int_{\Omega}  \eta^2  (1 + |Du|^2)^{\frac{\gamma+p(x)-2}{2}}|D^2u|^2 \, dx  \cr\cr 
& & + \left(\frac{2}{\sigma}+1\right) \int_{\Omega}  \eta^2   k^2(x)(1 + |Du|^2)^{\frac{\gamma+p(x)}{2}}\log^2(e+|Du|) \, dx  \cr\cr 
&&\quad +  \, \int_{\Omega}|D \eta|^2   (1 + |Du|^2)^{\frac{\gamma+p(x)}{2}}\,dx. \label{(4.10)}
\end{eqnarray}
Inserting estimates \eqref{I1} and \eqref{(4.10)} in \eqref{ristart}, we obtain
\begin{eqnarray*}
&&  \int_{\Omega} \eta^2  (1 + |Du|^2)^{\frac{\gamma}{2}}  \, \sum_{i,j,s, \alpha, \beta} F_{\xi_i^{\alpha} \xi_j^{\beta}}(x, Du) u_{x_j x_s}^{\alpha} u_{x_i x_s}^{\beta} \, dx \nonumber\\
&\le&  \frac{1}{2} \int_{\Omega} \eta^2  (1 + |Du|^2)^{\frac{\gamma}{2}}  \, \sum_{i,j,s, \alpha, \beta} F_{\xi_i^{\alpha} \xi_j^{\beta}}(x, Du) u_{x_j x_s}^{\alpha} u_{x_i x_s}^{\beta} \, dx
\cr\cr
&&\quad+\sigma\int_{\Omega}  \eta^2  (1 + |Du|^2)^{\frac{\gamma+p(x)-2}{2}}|D^2u|^2 \, dx  \cr\cr 
& &\quad + \left(\frac{2}{\sigma}+1\right) \,\int_{\Omega}  \eta^2   k^2(x)(1 + |Du|^2)^{\frac{\gamma+p(x)}{2}}\log^2(e+|Du|) \, dx  \cr\cr 
&&\quad + \, (2\ell+1) \int_{\Omega} |D \eta|^2  (1 + |Du|^2)^{\frac{\gamma+p(x)}{2}}\,dx.
\end{eqnarray*}
Reabsorbing the first integral in the right hand side by the left hand side, we get
\begin{eqnarray*}
&& \int_{\Omega} \eta^2  (1 + |Du|^2)^{\frac{\gamma}{2}}  \, \sum_{i,j,s, \alpha, \beta} F_{\xi_i^{\alpha} \xi_j^{\beta}}(x, Du) u_{x_j x_s}^{\alpha} u_{x_i x_s}^{\beta} \, dx \nonumber\\
&\le& 2\sigma\int_{\Omega}  \eta^2  (1 + |Du|^2)^{\frac{\gamma+p(x)-2}{2}}|D^2u|^2 \, dx  \cr\cr 
& &\quad + \, 2\left(\frac{2}{\sigma}+1\right) \int_{\Omega}  \eta^2   k^2(x)(1 + |Du|^2)^{\frac{\gamma+p(x)}{2}}\log^2(e+|Du|) \, dx  \cr\cr 
&&\quad + 2(2\ell+1) \int_{\Omega}   (1 + |Du|^2)^{\frac{\gamma+p(x)}{2}}|D \eta|^2\,dx
\end{eqnarray*}
and so, by the ellipticity assumption (H2), we deduce
\begin{eqnarray*}
&& \nu\int_{\Omega} \eta^2  (1 + |Du|^2)^{\frac{\gamma+p(x)-2}{2}}|D^2u|^2 \, dx \nonumber\\
&\le& 2\sigma \, \int_{\Omega}  \eta^2  (1 + |Du|^2)^{\frac{\gamma+p(x)-2}{2}}|D^2u|^2 \, dx  \cr\cr 
& &\quad + 2\left(\frac{2}{\sigma}+1\right) \int_{\Omega}  \eta^2   k^2(x)(1 + |Du|^2)^{\frac{\gamma+p(x)}{2}}\log^2(e+|Du|) \, dx  \cr\cr 
&&\quad + 2(2\ell+1) \int_{\Omega}  |D \eta|^2 (1 + |Du|^2)^{\frac{\gamma+p(x)}{2}}\,dx.
\end{eqnarray*}
Choosing $\sigma=\frac{\nu}{4}$, we can reabsorb the first integral in the right hand side by the left hand side {in order to} obtain
\begin{eqnarray}\label{4.10bis}
&& \int_{\Omega} \eta^2  (1 + |Du|^2)^{\frac{\gamma+p(x)-2}{2}}|D^2u|^2 \, dx \nonumber\\
&\le&  \, \bar C \, \int_{\Omega}  \eta^2   k^2(x)(1 + |Du|^2)^{\frac{\gamma+p(x)}{2}}\log^2(e+|Du|) \, dx  \cr\cr 
&&\quad + \, \bar C \, \int_{\Omega} |D \eta|^2  (1 + |Du|^2)^{\frac{\gamma+p(x)}{2}}\,dx,
\end{eqnarray}
with $\bar C=\frac{8}{\nu}\left(\frac{2}{\nu}+\ell+1\right).$
\\
One can easily check that
\begin{eqnarray*}
&& \left |D \left ( \eta \left (1 + |Du|^2\right )^{\frac{\gamma + p(x)}{4}} \right ) \right |^2\\
&\le& \, C \, (\gamma + p(x))^2 \, \eta^2 \, (1 + |Du|^2)^{\frac{\gamma + p(x) -2}{2}} |D^2 u|^2 \\
&&\qquad+ \eta^2 \, (1 + |Du|^2)^{\frac{\gamma + p(x)}{2}} k^2(x) \, \log^2 (1 + |Du|^2) + |D \eta|^2 \,  \left (1 + |Du|^2\right )^{\frac{\gamma + p(x)}{2}}\\
&\le& \, C \, (\gamma + p^+)^2 \, \eta^2 \, (1 + |Du|^2)^{\frac{\gamma + p(x) -2}{2}} |D^2 u|^2 \\
&&\qquad+ \eta^2 \, (1 + |Du|^2)^{\frac{\gamma + p(x)}{2}} k^2(x) \, \log^2 (1 + |Du|^2)  + |D \eta|^2 \,  \left (1 + |Du|^2\right )^{\frac{\gamma + p(x)}{2}} \, dx,
\end{eqnarray*}
where, in the last inequality, we used \eqref{p+p-} and $C$ is an absolute constant.
Now integrating both sides of the previous estimate in $\Omega$, using  \eqref{4.10bis} {and taking into account that $\log^2(1 + |Du|^2) \le \, C \, \log^2(e + |Du|),$} we obtain
\begin{eqnarray*}
&& \int_{\Omega} \left |D \left ( \eta \left (1 + |Du|^2\right )^{\frac{\gamma + p(x)}{4}} \right ) \right |^2 \, dx \\
&\le& \,C(\gamma + p^+)^2\int_{\Omega}  \eta^2   k^2(x)(1 + |Du|^2)^{\frac{\gamma+p(x)}{2}}\log^2(e+|Du|) \, dx
\\
&&\quad+ C \, (\gamma + p^+)^2 \,  \int_{\Omega} |D \eta|^2 (1 + |Du|^2)^{\frac{\gamma + p(x)}{2}},
\end{eqnarray*}
where we used that $(\gamma + p^+)^2\ge 1.$
Being $n>2$, by the Sobolev Embedding theorem, we deduce 
\begin{eqnarray*}
&& \left (\int_{\Omega} \left |\left ( \eta \left (1 + |Du|^2\right )^{\frac{\gamma + p(x)}{4}} \right ) \right |^{2^*} \, dx \right )^{\frac{2}{2^*}}\\
&\le& \, C \, (\gamma + p^+)^2  \int_{\Omega} \eta^2 \, (1 + |Du|^2)^{\frac{\gamma + p(x)}{2}} k^2(x) \, {\log^2 (e + |Du|)} \, dx \cr\cr
&&\quad+ \,C \, (\gamma + p^+)^2  \int_{\Omega} |D \eta|^2 \,  \left (1 + |Du|^2\right )^{\frac{\gamma + p(x)}{2}} \, dx,
\end{eqnarray*}
where $\displaystyle{\frac{2^*}{2} = \frac{n}{n-2}}.$ This yields
\begin{eqnarray}\label{semifinale}
\nonumber
&& \left (\int_{\Omega}  \eta^{2^*} \left (1 + |Du|^2\right )^{\frac{\gamma + p(x)}{2}\frac{2^*}{2}} \, dx \right )^{\frac{2}{2^*}}\\\nonumber
&\le& \, C \, (\gamma + p^+)^2  \int_{\Omega} \eta^2 \, (1 + |Du|^2)^{\frac{\gamma + p(x)}{2}} k^2(x) \, {\log^2 (e + |Du|)} \, dx \cr\cr
&&\quad+ \, C \, (\gamma + p^+)^2  \int_{\Omega} |D \eta|^2 \,  \left (1 + |Du|^2\right )^{\frac{\gamma + p(x)}{2}} \, dx\\
&=:& C \, (\gamma + p^+)^2\left(J_1+J_2\right).
\end{eqnarray}
For further needs, we record that 
in \eqref{semifinale} 
$C=C(n)\left(\frac{1}{\nu^2}+\frac{\ell+1}{\nu}\right).$
\\
%
We now split $J_1$ as follows
\begin{eqnarray}\label{J1}
\nonumber
J_1&= & \int_{\left \{ k {\log(e +|Du|)} \le \exp \left (\frac{1}{\varepsilon} \right ) - e \right \}} \eta^2 \, (1 + |Du|^2)^{\frac{\gamma + p(x)}{2}} k^2(x) \, \log^2 {(e +|Du|)} \, dx  \\\nonumber
&& + \int_{\left \{ k\log{(e +|Du|)} > \exp \left (\frac{1}{\varepsilon} \right ) - e \right \}} \eta^2 \, (1 + |Du|^2)^{\frac{\gamma + p(x)}{2}} k^2(x) \, \log^2 {(e +|Du|)} \, dx  \\\nonumber
&\le & \left(\exp \left (\frac{1}{\varepsilon} \right ) - e\right)^2\int_{\Omega} \eta^2 \, (1 + |Du|^2)^{\frac{\gamma + p(x)}{2}} dx  \\\nonumber
&& + \left(\int_{\Omega} \eta^2 \, (1 + |Du|^2)^{\frac{\gamma + p(x)}{2}\frac{n}{n-2}}  \, dx\right)^{\frac{n-2}{n}}\\
&&\qquad \cdot\left(\int_{\left \{ k {\log(e +|Du|)} > \exp \left (\frac{1}{\varepsilon} \right ) - e \right \}} \eta^2 \,  k^n(x) \, {\log^n(e +|Du|)} \, dx\right)^{\frac{2}{n}}
\end{eqnarray}
where $\varepsilon>0$ is a parameter that will be chosen in the sequel. We now have
\begin{eqnarray}
\nonumber
&&\int_{\left \{ k {\log(e +|Du|)} > \exp \left (\frac{1}{\varepsilon} \right ) - e \right \}} \eta^2 \,  k^n(x) \, \log^n(e +|Du|) \, dx
\\\nonumber
 &=& \int_{\left \{ k {\log(e +|Du|)} > \exp \left (\frac{1}{\varepsilon} \right ) - e \right \}} \eta^2 k^n(x) \frac{\log^{\theta} (e + k{\log(e +|Du|)}) }{\log^{\theta} (e + k {\log(e +|Du|)})} \log^n (e + |Du|) \, dx 
\\\nonumber
&\le & \varepsilon^{\theta} \int_{\Omega} \eta^2 k^n(x) \log^{\theta}(e + k {\log(e +|Du|)}) \log^n(e + |Du|) \, dx
\\\nonumber
&\le & (2\varepsilon)^{\theta} \int_{\Omega} \eta^2 k^n(x) [\log^{\theta}(e + k)+\log^{\theta}(e+ {\log(e +|Du|)})]  \log^n(e + |Du|) \, dx\\\nonumber
&\le& (2\varepsilon)^{\theta} \int_{\Omega} \eta^2 k^n(x)\log^{\theta}(e + k)  \log^n(e + |Du|) \, dx \\ \nonumber
&& + \, C \, \varepsilon^{\theta} \int_{\Omega} \eta^2 k^n  \log^{n+\theta}(e + |Du|) \, dx
\\
&=:& A_1+A_2,\label{intbrutto}
\end{eqnarray}
where $\theta$ is the exponent in assumption (H5) and the constant $C$ depends only on $\theta$.
\\
At this point, we use the elementary inequality stated in Lemma \ref{elem} with the choices 
\[
\delta = \tau = 1 \qquad \beta = n \qquad s = k^n(x) \log^{\theta}(e + k(x)) \qquad t = \log^n(e + |Du|)
\]
to deduce that
\begin{eqnarray*} 
A_1 &\le& (2 \varepsilon)^{\theta} \int_{\Omega} \eta^2 k^n(x) \log^{\theta}(e + k(x)) \log^n [e + k^n \log^{\theta}(e + k(x))] \, dx \\
&& + (2 \varepsilon)^{\theta} \int_{\Omega} \eta^2 (e+|Du|)\log(e + |Du|) \, dx \\
& =:& A_{11} + A_{12}.
\end{eqnarray*}
Since $\log^{\theta}(e + k^n(x)) > 1,$ we have that
\[
e + k^n(x) \log^{\theta} (e + k^n(x)) \le \, (e + k^n(x)) \log^{\theta} (e + k^n(x))
\]
and so, by elementary properties of the logarithm
\begin{eqnarray*}
&& \log^n [e + k^n(x) \log^{\theta} (e + k^n(x))] \\
&\le \, & \log^n [ (e + k^n(x)) \log^{\theta} (e + k^n(x))]
\\
&\le& C(n) \log^n(e + k^n(x)) + \log^n(e + \log^{\theta} (e + k^n(x))) \\
&\le & C(n) \log^n(e + k^n(x)) + {C(n, \theta)}\log^{n}(e + k^n(x))
\end{eqnarray*}
therefore
\begin{eqnarray*}
A_{11} 
\le  C\varepsilon^{\theta} \int_{\Omega} \eta^2 k^n(x) \log^{n + \theta}(e+k^n(x)) \, dx 
\end{eqnarray*}
and
\begin{eqnarray}\label{A1}
 A_1&\le& \, C \, \varepsilon^{\theta} \int_{\Omega} \eta^2 k^n(x) \log^{n + \theta}(e + k^n(x)) \, dx\cr\cr
 &&+C\varepsilon^{\theta} \int_{\Omega} \eta^2 (e+|Du|)\log(e + |Du|) \, dx.
\end{eqnarray}
On the other hand, using once more Lemma \ref{elem} with the choices
\[
\delta = \tau = 1 \qquad \beta = n+\theta \qquad s = k^n(x) \qquad t = \log^{n+\theta}(e + |Du|)
\] 
we obtain 
\begin{eqnarray}\label{A2}
A_2 &\le& C\varepsilon^\theta\int_{\Omega} \eta^2 k^n(x) \log^{n+\theta} (e + k^n(x)) \, dx\cr\cr
&&\quad + C\varepsilon^\theta\int_{\Omega} \eta^2 (e+|Du|)\log^{n+\theta} (e +|Du|) \, dx.
\end{eqnarray}
Inserting \eqref{A1} and \eqref{A2} in \eqref{intbrutto}, we obtain 
\begin{eqnarray}
\nonumber
&&\int_{\left \{ k {\log(e +|Du|)} > \exp \left (\frac{1}{\varepsilon} \right ) - e \right \}} \eta^2 \,  k^n(x) \, {\log^n(e +|Du|)} \, dx
\\\nonumber
&\le& \, C \, \varepsilon^{\theta} \int_{\Omega} \eta^2 k^n(x) \log^{n + \theta}(e + k^n(x)) \, dx + C \, \varepsilon^{\theta} \int_{\Omega} \eta^2 (e + |Du|) \log^{n+\theta}(e + |Du|) \, dx.
\end{eqnarray}
Therefore, recalling \eqref{J1}, we have 
\begin{eqnarray}\label{J1bis}
\nonumber
J_1
&\le & \left(\exp \left (\frac{1}{\varepsilon} \right ) - e\right)^{{2}} \int_{\Omega} \eta^2 \, (1 + |Du|^2)^{\frac{\gamma + p(x)}{2}} dx  \\
&& + \, C \, \varepsilon^{\frac{2\theta}{n}}\left(\int_{\Omega} \eta^2 \, (1 + |Du|^2)^{\frac{\gamma + p(x)}{2}\frac{n}{n-2}}  \, dx\right)^{\frac{n-2}{n}}\\\nonumber
&& \cdot \left(\int_{\Omega} \eta^2 k^n(x) \log^{n + \theta}(e + k^n(x)) \, dx+ \int_{\Omega} \eta^2 (e+|Du|)\log^{n+\theta}(e + |Du|) \, dx\right)^{\frac{2}{n}}
\end{eqnarray}
By virtue of (H5) and \eqref{embed}, recalling that $\eta\in C^\infty_0{(B_t)}$, with $B_t \Subset B_{R_0},$ we have that
\begin{eqnarray}\label{stimaDu}
&&\int_{\Omega} \eta^2 (e+|Du|)\log^{n+\theta}(e + |Du|) \, dx\le \int_{B_{R_0}} (e+|Du|)\log^{n+\theta}(e + |Du|) \, dx\cr\cr
&=&|B_{R_0}|\medint_{B_{R_0}} (e+|Du|)\log^{n+\theta}(e + |Du|) \, dx\cr\cr
&\le& C(n,\theta, p^-)|B_{R_0}|\left(\medint_{B_{R_0}} (e+|Du|)^{p^-}\,dx\right)^{\frac{1}{p^-}}\cr\cr
&\le&C(n,\theta, p^-)|B_{R_0}|\left(\medint_{B_{R_0}}(e+|Du|)^{p(x)}\,dx\right)^{\frac{1}{p^-}}\cr\cr
&\le& C(n,\theta,p^-, p^+, R_0)\left(1+\int_{B_{R_0}}|Du|^{p(x)}\,dx\right)^{\frac{1}{p^-}}\cr\cr
&\le& C(n,\theta,p^-, p^+, R_0)\left(1+M_1\right)^{\frac{1}{p^-}} =: \hat C\left(1+M_1\right)^{\frac{1}{p^-}},
\end{eqnarray}
 where, in the last inequality, we used  \eqref{M1}.
Setting 
\begin{equation}
\label{MR0}
M_{R_0}=||k||_{L^n\log^{n+\theta}L(B_{R_0})}+ \hat C \left(1+M_1\right)^{\frac{1}{p^-}},
\end{equation}
estimate \eqref{J1bis} becomes
\begin{eqnarray}\label{J1b}
J_1
&\le & \left(\exp \left (\frac{1}{\varepsilon} \right ) - e\right)^{{2}}\int_{\Omega} \eta^2 \, (1 + |Du|^2)^{\frac{\gamma + p(x)}{2}} dx  \\\nonumber
&& + \, C \, \varepsilon^{\frac{2\theta}{n}}M^{\frac{2}{n}}_{R_0}\left(\int_{\Omega} \eta^2 \, (1 + |Du|^2)^{\frac{\gamma + p(x)}{2}\frac{n}{n-2}}  \, dx\right)^{\frac{n-2}{n}}.
\end{eqnarray}
Inserting \eqref{J1b} in \eqref{semifinale}, we get
\begin{eqnarray}\label{prefinale}
&& \left (\int_{\Omega}  \eta^\frac{2n}{n-2} \left (1 + |Du|^2\right )^{\frac{\gamma + p(x)}{2}\frac{n}{n-2}} \, dx \right )^{\frac{n-2}{n}}\\\nonumber
&\le& \, C \, (\gamma + p^+)^2  \left(\exp \left (\frac{1}{\varepsilon} \right ) - e\right)^{{2}}\int_{\Omega} \eta^2 \, (1 + |Du|^2)^{\frac{\gamma + p(x)}{2}} dx  \\\nonumber
&& + \, C \, \varepsilon^{\frac{2\theta}{n}}(\gamma + p^+)^2 \, M^{\frac{2}{n}}_{R_0}\left(\int_{\Omega} \eta^2 \, (1 + |Du|^2)^{\frac{\gamma + p(x)}{2}\frac{n}{n-2}}  \, dx\right)^{\frac{n-2}{n}} \cr\cr
&&\quad+C \, (\gamma + p^+)^2  \int_{\Omega} |D \eta|^2 \,  \left (1 + |Du|^2\right )^{\frac{\gamma + p(x)}{2}} \, dx,
\end{eqnarray}
where $C=C(n,\theta)\cdot\bar C=C(n,\theta)\left(\frac{1}{\nu^2}+\frac{\ell+1}{\nu}\right)$.
Choose now $\varepsilon$ such that
$$C\varepsilon^{\frac{2\theta}{n}}(\gamma + p^+)^2M^{\frac{2}{n}}_{R_0}=\frac{1}{2} \,\,\Longleftrightarrow\,\, \varepsilon=\left(\frac{1}{2CM^{\frac{2}{n}}_{R_0}(\gamma + p^+)^2}\right)^{\frac{n}{2\theta}}\sim \left(\frac{1}{\Gamma(\gamma + p^+)^{\frac{n}{\theta}}}\right)$$
with $\Gamma=\Gamma(n,\theta,\nu, \ell,M_1, ||k||_{L^n\log^{n+\theta}L(B_{R_0})}, R_0, p^-)= C(n,\theta)M_{R_0}^{\frac{2}{n}}\left(\frac{1}{\nu^2}+\frac{\ell+1}{\nu}\right)$. With this choice \eqref{prefinale} and by the properties of the cut-off function $\eta$, becomes
\begin{eqnarray}\label{prefinale2}
\nonumber
&& \left (\int_{B_{{s}}}   \left (1 + |Du|^2\right )^{\frac{\gamma + p(x)}{2}\frac{2^*}{2}} \, dx \right )^{\frac{2}{2^*}}\\\nonumber
&\le& \,  \frac{1}{2}\left(\int_{B_{{t}}}  \, (1 + |Du|^2)^{\frac{\gamma + p(x) }{2}\frac{2^*}{2}}  \, dx\right)^{\frac{2}{2^*}}\cr\cr
&&+\frac{C}{{(t-s)^2}} \,  \left((\gamma + p^+)\exp \left( \Gamma(\gamma + p^+)^{\frac{n}{\theta}} \right) \right)^2\int_{B_{{t}}}  \, (1 + |Du|^2)^{\frac{\gamma + p(x)}{2}} dx,
\end{eqnarray}
since we may suppose without loss of generality that $t-s\le 1$ and we used that $\exp \left( \Gamma(\gamma + p^+)^{\frac{n}{\theta}} \right)\ge 1.$
\\
With the use of  Lemma \ref{holf}, we deduce that
\begin{eqnarray*}
&& \left (\int_{B_\rho}   \left (1 + |Du|^2\right )^{\frac{\gamma + p(x)}{2}\frac{2^*}{2}} \, dx \right )^{\frac{2}{2^*}}\\
&\le& \frac{C}{(R-\rho)^2} \,  \left((\gamma + p^+)\exp \left( \Gamma(\gamma + p^+)^{\frac{n}{\theta}} \right) \right)^2\int_{B_R}  \, (1 + |Du|^2)^{\frac{\gamma + p(x)}{2}} dx.
\end{eqnarray*}
This yields
\begin{eqnarray*}
&& \left (\int_{B_\rho}   \left (1 + |Du|^2\right )^{\frac{(\gamma + p^-)}{2} \frac{2^*}{2}} \, dx \right )^{\frac{2}{2^*}}\\
&\le& \frac{C}{(R-\rho)^2} \,  \left((\gamma + p^+)\exp \left( \Gamma(\gamma + p^+)^{\frac{n}{\theta}} \right) \right)^2\int_{B_R}  \, (1 + |Du|^2)^{\frac{\gamma + p^+}{2}} dx\\
&\le& \frac{C}{(R-\rho)^2} \,  \left(\exp \left( \Gamma(\gamma + p^+)^{\frac{n}{\theta}} \right) \right)^{2\left(1+\frac{\theta}{n}\right)}\int_{B_R}  \, (1 + |Du|^2)^{\frac{\gamma + p^+}{2}} \, dx,
\end{eqnarray*}
where we used that $t\le \exp(\Gamma t)$ and that $t=\left(t^{\frac{n}{\theta}}\right)^{\frac{\theta}{n}}\le \left(\exp(\Gamma t^{\frac{n}{\theta}})\right)^{\frac{\theta}{n}}$.
This in turn implies
\begin{eqnarray}\label{semifinito}
&& \left (\int_{B_\rho}   \left (1 + |Du|^2\right )^{\frac{(\gamma + p^-)}{2} \frac{2^*}{2}} \, dx \right )^{\frac{2}{2^*}}\\
&\le& \frac{C}{(R-\rho)^2} \,  \left(\exp \left (\Gamma(\gamma + p^-)^{\frac{n}{\theta}} + \Gamma(p^+ - p^- )^{\frac{n}{\theta}} \right ) \right)^{2\left(1+\frac{\theta}{n}\right)} \nonumber \\
&& \times \, \|(1 + |Du|^2)\|^{p^+ - p^-}_{L^{\infty}(B_R)} \int_{B_R}  \, (1 + |Du|^2)^{\frac{ \gamma+p^- }{2}} dx  \nonumber\\
&\le& \frac{C}{(R-\rho)^2} \,  \left(\exp \left (\Gamma(p^+ - p^-)^{\frac{n}{\theta}} + \Gamma (\gamma + p^-)^{\frac{n}{\theta}} \right ) \right)^{2\left(1+\frac{\theta}{n}\right)}\nonumber \\
&& \times \, \|(1 + |Du|^2)\|^{\omega(R_0)}_{L^{\infty}(B_R)} \int_{B_R}  \, (1 + |Du|^2)^{\frac{\gamma+p^- }{2}} dx  \nonumber\\
&\le& \frac{C}{(R-\rho)^2} \,  C \, \left(\exp \left (\Gamma(\gamma + p^-)^{\frac{n}{\theta}} \right ) \right)^{2\left(1+\frac{\theta}{n}\right)}  \nonumber\\
&& \times \, \|(1 + |Du|^2)\|^{\omega(R_0)}_{L^{\infty}(B_R)} \int_{B_R}  \, (1 + |Du|^2)^{\frac{\gamma + p^-}{2}} dx, \nonumber
\end{eqnarray}
where $C$ depends on $n,\theta, \ell, \nu, M_1, {||k||_{L^n\log^{n+\theta}L(B_{R_0})}}, R_0, p^-,p^+.$

\color{black}
Let us now define the decreasing sequence of radii $\rho_i$, $i \in \mathbb{N}$ by setting
\[
\rho_i := \rho + \frac{R - \rho}{2^i}
\]
and the increasing sequence of exponents
\[
p_i := p^- \left (\frac{2^*}{2} \right )^i.
\]
In this manner we can rewrite the  inequality  \eqref{semifinito} on every ball $B_{\rho_i}$ as follows
\begin{eqnarray*}
&& \left (\int_{B_{\rho_{i+1}}}   \left (1 + |Du|^2\right )^{\frac{p_{i+1}}{2}} \, dx \right )^{\frac{1}{p_{i+1}}}\\
&\le& \frac{C^{\frac{1}{p_i}}}{(\rho_i-\rho_{i+1})^{\frac{2}{p_i}}} \, 
\left(\exp \left (\Gamma(p_i)^{\frac{n}{\theta}} \right ) \right)^{\frac{\tilde\vartheta}{p_i}} \|(1 + |Du|^2)\|^{\frac{\omega(R_0)}{p_i}}_{L^{\infty}(B_R)} \left (\int_{B_{\rho_i}}  \, (1 + |Du|^2)^{\frac{p_i}{2}} dx \right )^{\frac{1}{p_i}},
\end{eqnarray*}
where {we used that  $\|(1 + |Du|^2)\|_{L^{\infty}(B_{\rho_i})}  \le \,  \|(1 + |Du|^2)\|_{L^{\infty}(B_R)}$ and where we set $\tilde\vartheta=2\left(1+\frac{\theta}{n}\right)$ in order to simplify the notation}.

By iterating the previous estimate, {for $i = 1, \dots, m,$} we obtain
\begin{eqnarray*}
&& \left (\int_{B_{\rho_{m+1}}}   \left (1 + |Du|^2\right )^{\frac{p_{m+1}}{2}} \, dx \right )^{\frac{1}{p_{m+1}}}\\
&\le& \prod_{i = 1}^m \frac{C^{\frac{1}{p_i}}}{(\rho_i-\rho_{i+1})^{\frac{2}{p_i}}} \, 
\left(\exp \left (\Gamma(p_i)^{\frac{n}{\theta}} \right ) \right)^{\frac{\tilde\vartheta}{p_i}} \|(1 + |Du|^2)\|^{\frac{\omega(R_0)}{p_i}}_{L^{\infty}(B_R)} \left (\int_{B_{R}}  \, (1 + |Du|^2)^{\frac{p^-}{2}} dx \right )^{\frac{1}{p^-}}\\
&\le& \prod_{i = 1}^m \frac{C^{\frac{1}{p_i}} \, 4^{\frac{i+1}{p_i}}}{(R-\rho)^{\frac{2}{p_i}}} \, 
\left(\exp \left (\Gamma(p_i)^{\frac{n}{\theta}} \right ) \right)^{\frac{\tilde\vartheta}{p_i}} \|(1 + |Du|^2)\|^{\frac{\omega(R_0)}{p_i}}_{L^{\infty}(B_R)} \left (\int_{B_{R}}  \, (1 + |Du|^2)^{\frac{p^-}{2}} dx \right )^{\frac{1}{p^-}}.
\end{eqnarray*}
At this point we observe that 
\begin{eqnarray*}
&& \prod_{i = 1}^m \frac{C^{\frac{1}{p_i}} \, 4^{\frac{i+1}{p_i}}}{(R-\rho)^{\frac{2}{p_i}}} \, 
\left(\exp \left (\Gamma(p_i)^{\frac{n}{\theta}} \right ) \right)^{\frac{\tilde\vartheta}{p_i}} \|(1 + |Du|^2)\|^{\frac{\omega(R_0)}{p_i}}_{L^{\infty}(B_R)} \\
&= &\prod_{i = 1}^m \frac{C^{\frac{1}{p_i}} \, 4^{\frac{i+1}{p_i}}}{(R-\rho)^{\frac{2}{p_i}}} \,  \prod_{i = 1}^m
\left(\exp \left (\Gamma(p_i)^{\frac{n}{\theta}} \right ) \right)^{\frac{\tilde\vartheta}{p_i}}  \prod_{i = 1}^m \|(1 + |Du|^2)\|^{\frac{\omega(R_0)}{p_i}}_{L^{\infty}(B_R)}.
\end{eqnarray*}
First of all
\begin{eqnarray*}
\prod_{i = 1}^m \frac{C^{\frac{1}{p_i}} \, 4^{\frac{i+1}{p_i}}}{(R-\rho)^{\frac{2}{p_i}}} &=& \exp \left (\log \frac{C}{(R - \rho)^2} \left(\sum_{i=0}^m \frac{1}{p_i} +  \sum_{i = 0}^m \frac{i + 1}{p_i} \right)\right )\\
&\le \, & \exp \left ( c \log \frac{C}{(R - \rho)^2} \, \left [\sum_{i=0}^m \frac{1}{p_i} + \sum_{i=0}^m \frac{i + 1}{p_i} \right ]\right )
\\
&\le& \exp \left ( c(\sigma) \log \frac{C}{(R - \rho)^2}  \sum_{i=0}^{+ \infty}  \left (\frac{2}{2^*} \right )^{i \sigma}\right) \\
&\le&  \exp \left ( c(n,\sigma) \log \frac{C}{(R - \rho)^2}  \right) \le  \left( \frac{C}{(R - \rho)^2}  \right)^{c(n,\sigma)},
\end{eqnarray*}
for some $\sigma \in (0,1)$ (for simplicity we fix $\sigma = \frac{1}{2}$). On the other hand
\begin{eqnarray*}
 \prod_{i = 1}^m
\left(\exp \left (\Gamma(p_i)^{\frac{n}{\theta}} \right ) \right)^{\frac{\tilde\vartheta}{p_i}} &=& \exp \left (\sum_{i = 0}^m \frac{\tilde\vartheta}{p_i}\log \left(\exp \left (\Gamma(p_i)^{\frac{n}{\theta}} \right ) \right) \right )\\
&=&  \exp \left (\tilde\vartheta\Gamma \sum_{i=0}^m (p_i)^{\frac{n}{\theta} - 1}\right )\\
&=&  \exp \left (\tilde\vartheta\Gamma \sum_{i=0}^m \left (\frac{1}{p^-} \right)^{\frac{n}{\theta} - 1} \left (\frac{2^*}{2}\right)^{\frac{n}{\theta} - 1}\right )\\
\end{eqnarray*}
and this is a convergent series provided that $n < \theta$.
Finally 
\begin{eqnarray*}
\prod_{i = 1}^m \|(1 + |Du|^2)\|^{\frac{\omega(R_0)}{p_i}}_{L^{\infty}(B_R)} &=& \exp \left (\sum_{i=0}^m \frac{\omega(R_0)}{p_i} \log \|(1 + |Du|^2)\|_{L^{\infty}(B_R)}\right ) \\
&=& \exp \left (\omega(R_0) \log \|(1 + |Du|^2)\|_{L^{\infty}(B_R)} \sum_{i=0}^m \frac{1}{p_i} \right ) \\
&\le & \exp \left (\omega(R_0) \log \|(1 + |Du|^2)\|_{L^{\infty}(B_R)} \frac{1}{p^-}\sum_{i=0}^{+\infty} \left (\frac{2}{2^*} \right )\right )\\
&=& \exp \left ( \omega(R_0) \log \|(1 + |Du|^2)\|_{L^{\infty}(B_R)} \frac{\alpha}{ p^-} \right ) \\
&=& \|(1 + |Du|^2)\|_{L^{\infty}(B_R)}^{\frac{\omega(R_0) \, \alpha}{ p^-}},
\end{eqnarray*}
where $\displaystyle \alpha = \frac{n}{2}.$
Therefore we end up with
\begin{eqnarray*}
&& \left (\int_{B_{\rho_{m+1}}}   \left(1 + |Du|^2\right)^{\frac{p_{m+1}}{2}} \, dx \right)^{\frac{1}{p_{m+1}}}\\
&\le & \hat C\|(1 + |Du|^2)\|^{\frac{\omega(R_0) \, \alpha}{ p^-}}_{L^{\infty}(B_R)} \left (\int_{B_{R}}  \, (1 + |Du|^2)^{\frac{p^-}{2}} dx \right) ^{\frac{1}{p^-}}
\end{eqnarray*}
for a suitable constant $\hat{C}$ depending on {$n,\theta,\nu, \ell,M_1, ||k||, R_0, p^-, p^+$}. Now letting $m \rightarrow + \infty$ we have
\[
||1+Du||_{L^\infty (B_\rho)}\le \hat{C}||1+Du||_{L^\infty(B_R)}^{\frac{\omega(R_0) \, \alpha}{ p^-}} \left (\int_{B_{R}}  \, (1 + |Du|^2)^{\frac{p^-}{2}} dx \right) ^{\frac{1}{p^-}}
\]
for every $\frac{R_0}{2}<\rho<R<R_0$. The continuity of $p(x)$ allows us to choose  $R_0$ such that $\omega(R_0)<\frac{1}{2}\frac{p^-}{\alpha}$ and so previous estimate becomes
\[
||1+Du||_{L^\infty (B_\rho)}\le \hat{C}||1+Du||_{L^\infty(B_R)}^{\frac{1}{ 2}} \left (\int_{B_{R}}  \, (1 + |Du|^2)^{\frac{p^-}{2}} dx \right) ^{\frac{1}{p^-}}
\]
and so, by Young's inequality, we get
\[
||1+Du||_{L^\infty (B_\rho)}\le \frac{1}{2}||1+Du||_{L^\infty(B_R)}+ \hat{C}\left (\int_{B_{R}}  \, (1 + |Du|^2)^{\frac{p^-}{2}} dx \right) ^{\frac{2}{p^-}}.
\]
The iteration Lemma \ref{holf} yields
\begin{equation}\label{apriori}
||1+Du||_{L^\infty \left (B_{\frac{R_0}{2}}\right )}\le  C\left (\int_{B_{R_0}}  \, (1 + |Du|^2)^{\frac{p^-}{2}} dx \right) ^{\frac{2}{p^-}},
\end{equation}
with a constant $C$ proportional to $\frac{1}{\nu^2}+\frac{\ell+1}{\nu}$, i.e. the conclusion.

\section{The approximation}

\label{quattro}

In this subsection we give only a few hints for the approximation procedure that comes along quite standard arguments (in view of the assumption \eqref{modulo}).

For every $x_0\in \Omega$ and every $\delta_0$ there exists a  ball $B_{R_0}(x_0)\Subset\Omega$ such  that, setting
\begin{equation}
\label{exp-p}
p=\inf\{p(x):\,\, x\in B_{R_0}\}>1
\end{equation}
\begin{equation}
\label{exp-q}
q=\sup\{p(x):\,\, x\in B_{R_0}\}
\end{equation}
one has
$$\frac{q}{p}<1+\delta_0.$$
Indeed
$$ q-p\le \omega(R_0)$$
and by the continuity assumption of $p(x)$ we have that for every $\delta>0$ there exists $R_\delta$ such that $\omega(R)<\delta$, for every $R<R_\delta$.
Therefore
$$\frac{q}{p}=\frac{q-p}{p}+1<\delta+1.$$
From \eqref{p+p-}, we deduce that
$$C_1 (1+|\xi|^2)^{\frac{p}{2}}\le F(x,\xi)\le \, C_2 (1+|\xi |^2)^{\frac{q}{2}}\eqno{\rm (H1')}$$
for a.e. $x\in B_{R_0}$, all $\xi\in \M$ and $0 < C_1 < C_2.$

 For a smooth mollifier $\varphi\in C^\infty_0(B_1(0))$ and for $K\in\mathbb{N} > 0$ let us introduce 
$$g_\varepsilon(x,\xi)=\int_{B_{1}(0)} F(x+\varepsilon y, \xi)\varphi(y)\,dy$$
and, following \cite{cupguimas}, the sequence of functionals
$$F_{\varepsilon, K}(x,\xi)=
\left\{
\begin{array}{lll} 
\!\!\! & g_\varepsilon(x,\xi) \qquad & |\xi| \le \, K\\[2mm]
\!\!\! & g_\varepsilon(x,K)+D_\xi g_\varepsilon(x,K)(|\xi|-K)\qquad & |\xi| \ge \, K
\end{array}
\right.$$
that are smooth with respect to the $x$ variable. 

The definition of $F_{\varepsilon,K}$, by virtue of the  continuity of $p(\cdot)$, entails
$$C_1(\varepsilon,K) (1+|\xi|^2)^{\frac{p(x)}{2}}\le F_{\varepsilon,K}(x,\xi)\le  C_2(\varepsilon, K) (1+|\xi |^2)^{\frac{p(x)}{2}}\eqno{\rm (A1)}$$
for a.e. $x\in B_{R_0}$, all $|\xi| \le \, K$, where 
\begin{eqnarray*}
&& C_1(\varepsilon, K) \sim C_1(1+K)^{-\log^{- \frac{\theta + 1}{n}} \left (e + \frac{1}{\varepsilon} \right)}
\\
&& C_2(\varepsilon, K) \sim C_2(1+K)^{\log^{- \frac{\theta + 1}{n}} \left (e + \frac{1}{\varepsilon} \right)}.
\end{eqnarray*}
For further needs, we observe that  
that $C_i(\varepsilon, K)  \rightarrow 1$ as $\varepsilon \rightarrow 0^+,$ for $i = 1,2$ and each fixed $k.$

\noindent Moreover assumptions (H2)--(H4) yield 
$$\langle {D_{\xi \xi}F_{\varepsilon,K}} (x,\xi)\eta,\eta\rangle \ge \nu \, C_1(\varepsilon, K)   (1 + |\xi |^2)^{\frac{p(x)-2}{2}}|\eta|^2 \eqno{\rm (A2)}$$
$$|D_{\xi \xi}F_{\varepsilon,K}(x,\xi)| \le  \ell \, C_2(\varepsilon, K)  (1 + |\xi |^2)^{\frac{p(x)-2}{2}}
\eqno{\rm (A3)}$$
$$\vert D_xD_\xi F_{\varepsilon,K}(x,\xi)\vert \le C_2(\varepsilon, K)k_\varepsilon (1+|\xi |^2)^{\frac{p(x)-1}{2}}\log(e+|\xi|^2)\eqno{\rm (A4)}$$
where $k_\varepsilon$ is the usual mollification of $k$.
\\
Obviously the use of \eqref{p+p-} in (A1)--(A4) yields the corresponding inequalities where from below we replace $p(x)$ with $p$ and from above with $q$, where $p$ and $q$ are defined in \eqref{exp-p} and \eqref{exp-q} respectively. 
\\
Now for $u\in W^{1,p(x)}(\Omega)$ a minimizer of the functional $\mathfrak{F}$, let $v_{\varepsilon,K}\in u+W^{1,p(x)}_0(B_{R_0})$ be the unique solution  to the problem
$$\min\left\{\mathfrak{F}_{\varepsilon,K} (v,B_{R_0}):\,\,\, v\in u+W^{1,p(x)}_0(B_{R_0})\right\},$$
where $$\mathfrak{F}_{\varepsilon,K} (v,B_{R_0})=\int_{B_{R_0}}F_{\varepsilon,K}(x, Dv)\,dx.$$
By virtue of (A2)--(A4) with $p$ and $q$ instead of $p(x),$  choosing $0<\delta_0<\frac{1}{n}$, we are legitimate to apply the result in \cite{EMM1} to obtain that $v_{\varepsilon,K}\in W^{1,\infty}_{\rm loc}(B_{R_0})$ and that in \cite{CMMP} to obtain that $v_{\varepsilon,K}\in W^{2,2}_{\rm loc}(B_{R_0}),$  for every $\varepsilon, K >0$.

The left inequality in (A1) implies that
\begin{eqnarray}
C_1(\varepsilon, K) \int_{B_{R_0}}(1+|Dv_{\varepsilon,K}|^2)^{\frac{p}{2}}\, dx&\le& \int_{B_{R_0}}F_{\varepsilon,K} (x,Dv_{\varepsilon,K})\,dx \label{normaLp} \\[2mm]
&\le& \int_{B_{R_0}}F_{\varepsilon,K} (x,Du)\,dx \nonumber \\[2mm]
&\le& C_2(\varepsilon, K) \int_{B_{R_0}}(1+|Du|^2)^{\frac{p(x)}{2}}\, dx  \nonumber
\end{eqnarray}
and so, using the explicit expression of $C_i(\varepsilon,K)$,
\[
\int_{B_{R_0}}(1+|Dv_{\varepsilon,K}|^2)^{\frac{p}{2}}\, dx \lesssim \,  \frac{C_2}{C_1}(1+K)^{2\log^{- \frac{\theta + 1}{n}} \left (e + \frac{1}{\varepsilon} \right)} \int_{B_{R_0}}(1+|Du|^2)^{\frac{p(x)}{2}}\, dx 
\]
that, taking the limit as $\varepsilon\to 0$, passing possibly to a not relabelled subsequence, in turn entails
\[
\int_{B_{R_0}}(1+|Dv_{\varepsilon,K}|^2)^{\frac{p}{2}}\, dx \lesssim \,  \widetilde{C}\, \int_{B_{R_0}}(1+|Du|^2)^{\frac{p(x)}{2}}\, dx 
\]
with $\widetilde{C}$ independent of $\varepsilon, K$. 
\\
Therefore $\big(v_{\varepsilon,K}\big)_\varepsilon$ is a bounded sequence in $W^{1,p}(B_{R_0})$ and so there exists $v_K\in W^{1,p}(B_{R_0})$ such that \begin{equation}\label{weakLp}
  v_{\varepsilon,K} \rightharpoonup v_K \qquad \text{in}\,\, W^{1,p}(B_{R_0}).  
\end{equation}
The lower semicontinuity of the $L^p$ norm implies
\[
\int_{B_{R_0}}(1+|Dv_{K}|^2)^{\frac{p}{2}}\, dx\le \liminf_{\varepsilon\to 0}\int_{B_{R_0}}(1+|Dv_{\varepsilon,K}|^2)^{\frac{p}{2}}\, dx \le \,  \widetilde{C}\, \int_{B_{R_0}}(1+|Du|^2)^{\frac{p(x)}{2}}\, dx 
\]
and so 
there exists $v\in W^{1,p}(B_{R_0})$ such that \begin{equation}\label{weakLpbis}
  v_{K} \rightharpoonup v \qquad \text{in}\,\, W^{1,p}(B_{R_0}).  
\end{equation}
On the other hand, since $v_{\varepsilon,K}\in W^{1,\infty}_{\rm loc}(B_{R_0}) \cap W^{2,2}_{\rm loc}(B_{R_0})$, we can use our a priori estimate at \eqref{apriori} to deduce that
\begin{eqnarray}\label{apriori2}
||1+Dv_{\varepsilon,K}||_{L^\infty \left (B_{\frac{R_0}{2}}\right )}&\le & C\left (\int_{B_{R_0}}  \, (1 + |Dv_{\varepsilon,K}|^2)^{\frac{p}{2}} dx \right) ^{\frac{2}{p}}\cr\cr
&\le&  C\left (\int_{B_{R_0}}  \, (1 + |Du|^2)^{\frac{p(x)}{2}} dx \right) ^{\frac{2}{p}},
\end{eqnarray}
with $C$ proportional to $C_1^{-1}(\varepsilon, K), C_2(\varepsilon, K)$ and where in the last inequality we used \eqref{normaLp}. Taking the limit  as  $\varepsilon\to 0$, \eqref{apriori2} holds with some constant $\overline{C}$ independent of $\varepsilon, K.$
\\
Therefore, passing to the limit as $K\to +\infty$, have that
\begin{equation}\label{weakLpbi}
  v_{\varepsilon,K} \stackrel{*}\rightharpoonup v \qquad \text{weakly * in}\,\, W^{1,\infty}_{\mathrm{loc}}(B_{R_0})  
\end{equation}
and, by the lower semicontinuity of the norm, also
\begin{eqnarray}\label{apriori2bis}
||1+Dv||_{L^\infty \left(B_{\frac{R_0}{2}}\right )}
&\le&  C\left (\int_{B_{R_0}}  \, (1 + |Du|^2)^{\frac{p(x)}{2}} dx \right) ^{\frac{2}{p}}.
\end{eqnarray}
At this point, the conclusion of the proof, that is showing that $v=u$ a.e. in $B_{R_0}$, follows by standard arguments that can be found for instance in \cite{EMM1}.

\vspace{10mm}

\noindent {\bf Aknowledgments.}
M. Eleuteri and A. Passarelli di Napoli have been partially supported by the Gruppo
Nazionale per l’Analisi Matematica, la Probabilità e le loro Applicazioni (GNAMPA) of the Istituto
Nazionale di Alta Matematica (INdAM). Moreover M. Eleuteri has been partially supported by PRIN 2020 ``Mathematics for industry 4.0 (Math4I4)'' (coordinator P. Ciarletta) while
A. Passarelli di Napoli has been partially supported by Università
degli Studi di Napoli Federico II through the Project FRA ( 000022-75-2021-FRA-PASSARELLI)

\end{document}